\documentclass[preprint,12pt,3p]{elsarticle}

\usepackage{amssymb}
\usepackage{amsthm}

\usepackage{amsmath}
\usepackage{mathtools}
\usepackage{amsfonts}
\usepackage{mathrsfs}
\usepackage{verbatim}
\usepackage{etoolbox}
\usepackage{hyperref}
\usepackage[none]{hyphenat}

\hypersetup{
    colorlinks=true,
    linkcolor=blue,
    filecolor=magenta,      
    urlcolor=cyan,
}

\newtheorem{thm}{Theorem}[section]
     
\newtheorem{prop}[thm]{Proposition}
\newtheorem{conj}[thm]{Conjecture}
\newtheorem{rmk}[thm]{Remark}

\numberwithin{equation}{section}

\newcommand{\N}{\mathbb{N}}
\newcommand{\Z}{\mathbb{Z}}
\newcommand{\R}{\mathbb{R}}

\journal{arXiv}

\begin{document}

\begin{frontmatter}

\title{New conjecture related to a conjecture of McIntosh \tnoteref{label1}}
\tnotetext[label1]{This research did not receive any specific grant from funding agencies in the public, commercial, or not-for-profit sectors.}

\author{Saud Hussein}
\address{Central Washington University, 400 East University Way, Ellensburg, WA 98926}

\ead{saudhussein23@gmail.com}

\begin{abstract}
We introduce a new conjecture on products of two distinct primes that would provide a partial answer to a conjecture of McIntosh \cite{McIntosh}. Also, $\binom{2p-1}{p-1}-1$ is written in terms of a polynomial in prime $p$ over the integers and we discuss one way this form may be useful.
\end{abstract}


\end{frontmatter}


\section{Introduction}

In 1771, Lagrange \cite{Lagrange} gave the first proof of an interesting property of the prime numbers we now call Wilson's theorem. The converse of this theorem also holds.

\begin{thm} [Wilson's theorem and its converse, theorem 5.4 \cite{Burton}] \label{wil}
\[(p-1)! \equiv -1 \pmod{p} \iff p \text{ is prime}.\]
\end{thm}

\begin{rmk}
For any $n \in \N$, \begin{align} \frac{(2n-1)!}{n!} = (n+1)(n+2)\cdots(n+n-1) \equiv (n-1)! \pmod{n}, \label{equ} \end{align} so theorem \ref{wil} may also be stated as \[\frac{(2p-1)!}{p!} \equiv -1 \pmod{p} \iff p \text{ is prime}.\]
\end{rmk}
\noindent
A prime $p$ satisfying the congruence \[(p-1)! \equiv -1 \pmod{p^2}\] is called a \textit{Wilson prime}. Infinitely many such primes are conjectured to exist but only $5$, $13$, and $563$ have been identified so far. We also think there are no integers $n$ such that \[(n-1)! \equiv -1 \pmod{n^3}.\]

For $n \in \N$ denote $w_n = \binom{2n-1}{n-1} = \frac{1}{2}\binom{2n}{n}$. Let $n=p$ be a prime in congruence \eqref{equ}. Since $p$ is relatively prime to $(p-1)!$, dividing by $(p-1)!$, \[w_p \equiv 1 \pmod{p}.\]
The congruence also holds for squares of odd primes, cubes of primes $\geq 5$, and for some products of distinct primes (see section \ref{sec3}). In 1819, Babbage \cite{Babbage} further showed \[w_p \equiv 1 \pmod{p^2}\] for primes $p \geq 3$. See conjecture \ref{mci} for other solutions. Finally, in 1862 Wolstenholme \cite{Wolstenholme} improved on Babbage's result.

\begin{thm} [Wolstenholme's theorem] If $p \geq 5$ is prime, then \[w_p \equiv 1 \pmod{p^3}.\]
\end{thm}
\noindent
James P. Jones conjectures no other solutions exist.

\begin{conj} [Jones' conjecture]
\[w_p \equiv 1 \pmod{p^3} \iff p \geq 5 \text{ is prime}.\]
\end{conj}
\noindent
Jones' conjecture is true for even integers, powers of primes $\leq$ $10^9$ (\cite{Trevisan}, \cite{Helou}), and based on computations for integers $\leq$ $10^9$. McIntosh \cite{McIntosh} gives probabilistic evidence supporting the conjecture (see conjecture \ref{mci}). A prime $p$ satisfying the congruence \[w_p \equiv 1 \pmod{p^4}\] is called a \textit{Wolstenholme prime}. McIntosh conjectures infinitely many such primes exist but only $16843$ and $2124679$ have been found so far. McIntosh also conjectures there are no integers $n$ such that \[w_n \equiv 1 \pmod{n^5}.\]

A summary of results known and conjectured related to Wilson's and Wolstenholme's theorems shows the similarity between the two notions.

\begin{center}
\begin{tabular}{ |c||c|c| } 
 \hline
 & \textbf{Wilson} & \textbf{Wolstenholme} \\
 m & $\frac{(2n-1)!}{n!} \equiv -1 \pmod{n^m}$ & $\binom{2n-1}{n-1} \equiv 1 \pmod{n^{2+m}}$ \\ 
 \hline
 \hline
 1 & if and only if $n$ is prime & if $n \geq 5$ is prime, conjectured only if \\ [0.5ex]
 \hline
 2 & $n = 5, 13, 563\dots$ conj $\infty$ many & $n = 16843, 2124679\dots$ conj $\infty$ many \\ [0.5ex]
 \hline
 3 & conjectured none & conjectured none \\ [0.5ex]
 \hline
\end{tabular}
\end{center}
\medskip
In section \ref{sec3}, we introduce a new conjecture that would provide a partial answer to a conjecture of McIntosh \cite{McIntosh}. In the last section, $(w_p-1)/p^3$ is written in terms of a polynomial in prime $p$ over the integers and we discuss one way this form may be useful. See \cite{Hussein} for a simple proof of Wilson's and Wolstenholme's theorems.

\section{Product of two distinct primes} \label{sec3}

From Bertrand's postulate and \cite[proposition 5, part 4]{Helou}, Jones' conjecture holds for products of two consecutive odd primes. For pairs of primes in general, the following is an equivalent criteria to the product satisfying Jones' conjecture.

\begin{prop} [corollary 4, proposition 4 \cite{Helou}]
Let $p$ and $q$ be distinct primes $\geq 5$. Then \[w_{pq} \equiv 1 \pmod{(pq)^3} \iff w_p \equiv 1 \pmod{q^3} \text{ and }  w_q \equiv 1 \pmod{p^3}.\]
\end{prop}
\noindent
A direct proof is outlined in \cite{Helou} that avoids the use of \cite[proposition 4]{Helou}. We fill in the details below. In \cite{McIntosh} for $n\in\N$, the \textit{modified binomial coefficient} is defined to be \[w_n'=\binom{2n-1}{n-1}'= \prod_{\substack{k=1 \\ (k,n)=1}}^n\frac{2n-k}{k}.\]
Notice $w_1'=1$. From \cite{Brinkmann} we also have the relation \begin{align} \binom{2n-1}{n-1} = \prod_{d|n} \binom{2d-1}{d-1}', \label{rel} \end{align} so $w_p=w_p'$ for any prime $p$.

\begin{proof}
By a generalization of Wolstenholme's theorem \cite[Theorem 1]{McIntosh}, \[w_{pq}' \equiv 1 \pmod{(pq)^3}.\]
Using relation \eqref{rel}, \[w_{pq} = w_1'w_p'w_q'w_{pq}' = w_pw_qw_{pq}',\] so \[w_{pq} \equiv w_pw_q \pmod{(pq)^3}.\]
We also have from Wolstenholme's theorem, \[w_p \equiv 1 \pmod{p^3} \quad \text{and} \quad w_q \equiv 1 \pmod{q^3},\] so by the previous congruence, \begin{align} w_{pq} \equiv w_q \pmod{p^3} \quad \text{and} \quad w_{pq} \equiv w_p \pmod{q^3}. \label{and} \end{align}

Now assume $w_{pq} \equiv 1 \pmod{(pq)^3}$. Then by \eqref{and}, \[w_p \equiv 1 \pmod{q^3} \quad \text{and} \quad w_q \equiv 1 \pmod{p^3}.\]
Conversely, assume \[w_p \equiv 1 \pmod{q^3} \quad \text{and} \quad w_q \equiv 1 \pmod{p^3}.\]
Then from \eqref{and}, \[w_{pq} \equiv 1 \pmod{p^3} \quad \text{and} \quad w_{pq} \equiv 1 \pmod{q^3}.\]
Since $p$ and $q$ are distinct primes, $w_{pq} \equiv 1 \pmod{(pq)^3}$ and the proof is complete.
\end{proof}

\begin{rmk} \label{rmk2}
Starting with \eqref{and} we can also show \[w_{pq} \equiv 1 \pmod{pq} \iff w_p \equiv 1 \pmod{q} \text{ and }  w_q \equiv 1 \pmod{p}\] and \[w_{pq} \equiv 1 \pmod{(pq)^2} \iff w_p \equiv 1 \pmod{q^2} \text{ and }  w_q \equiv 1 \pmod{p^2}.\]
\end{rmk}
\noindent
The only known examples of distinct primes $p$ and $q$ such that \[w_{pq} \equiv 1 \pmod{pq}\] are $(29,937)$, $(787, 2543)$, and $(69239, 231433)$. In \cite{Helou} we ask if there are actually infinitely many examples of such pairs of primes.

Let $n=p^2$ with $p$ a Wolstenholme prime. By the generalization of Wolstenholme's theorem \cite[theorem 1]{McIntosh}, \[w_n' \equiv 1 \pmod{n^2}.\]
Using relation \eqref{rel}, \[w_n = w_{p^2} = w_1'w_p'w_{p^2}' = w_pw_{p^2}' = w_pw_n',\] so \[w_n \equiv w_p \pmod{n^2}.\]
Since $p$ is a Wolstenholme prime, \[w_p \equiv 1 \pmod{(p^2)^2},\] so \[w_n \equiv 1 \pmod{n^2}.\]
Therefore McIntosh conjectures the following:

\begin{conj}[McIntosh's conjecture] \label{mci}
 \begin{align*} w_n \equiv 1 \pmod{n^2} \iff& n \text{ is odd prime or}\\
 &n=p^2 \quad \text{with } p \text{ a Wolstenholme prime.}
 \end{align*}
\end{conj}
\noindent
If McIntosh's conjecture is true, the only remaining composites that could violate Jones' conjecture are the Wolstenholme primes squared. By the definition of $w_p$ for prime $p\geq 5$, for any prime $q \geq \sqrt{2p-1}$, $q \not= p$, such that $\frac{2p-1}{n+1} < q \leq \frac{2p-1}{n}$, $n\in \N$, we have the following: \begin{align*} & n \text{ is odd} \implies q\mid w_p,\\
& n \text{ is even} \implies q \not \vert \,\, w_p.
\end{align*}
Hence most of the large prime factors of $w_p-1$ reside in one of the even $n$ inequalities above and calculations show the primes $q < p \leq 10^5$ such that $q^2\mid (w_p-1)$ are more than $100$ times smaller than $p$. Therefore it is increasingly unlikely for $q^2 \mid (w_p-1)$ for a prime $q > 2p$ the larger $p$ is, so we conjecture the following:

\begin{conj}[New conjecture] \label{conj}
For all but at most finitely many pairs of distinct primes $p$ and $q$, \[w_p \equiv 1 \pmod{q^2} \implies q < p.\]
\end{conj}
\noindent
The new conjecture implies \[w_p \not\equiv 1 \pmod{q^2} \quad \text{or} \quad w_q \not\equiv 1 \pmod{p^2}\] for all but at most finitely many pairs of distinct primes $p$ and $q$, so by Remark \ref{rmk2}, \[w_{pq} \not\equiv 1 \pmod{(pq)^2}.\]

\noindent
However, a prime pair $(p,q)$ exception in conjecture \ref{conj} only means $w_{pq} \equiv 1 \pmod{(pq)^2}$ is possible and an effective proof of conjecture \ref{conj} would allow us to check the finitely many exceptions. Hence the new conjecture would provide a partial answer to McIntosh's conjecture.

\section{Wolstenholme polynomials} \label{sec4}

For a prime $p \geq 5$, \begin{align} w_p = \binom{2p-1}{p-1} = \prod_{k=1}^{p-1}\left(\frac{p}{k}+1\right) = \sum_{k=0}^{p-1} p^k S_{p-1,k} \label{form}
\end{align}
where $S_{n,k}$ is the $k$th \textit{elementary symmetric polynomial} in $n$ variables with values $\{1,\frac{1}{2},\frac{1}{3},\dots,\frac{1}{n}\}$. Note $S_{n,0} = 1$, $S_{n,n} = 1/n!$, and $S_{n,k} = 0$ for $k > n$. By \cite[theorem 3]{Bayat} and Newton's identities relating power sums and elementary symmetric polynomials, \begin{align} \label{fra} S_{p-1,k} &\equiv 0 \pmod{p}, \quad \text{even } k \leq p-3\\
&\equiv 0 \pmod{p^2}, \quad \text{odd } k \leq p-4. \nonumber
\end{align}
The \textit{fractional congruence} notation \eqref{fra} means $p$, $p^2$ divides the numerator of the rational $S_{p-1,k}$ but not the denominator since the denominator is a product of integers $\leq p-1$. More generally, the highest power of $p$ dividing the numerator of $S_{p-1,k}$, $k=1,3,\dots,p-2$, is one higher than for $S_{p-1,k+1}$.

\begin{rmk}
Since $S_{p-1,p-1} = 1/(p-1)!$, $S_{p-1,p-2} \equiv 0 \pmod{p}$. Also by \eqref{form}, $S_{p-1,1} \equiv 0 \pmod{p^2}$ is equivalent to Wolstenholme's theorem and $p$ is a Wolstenholme prime if and only if $S_{p-1,1} \equiv 0 \pmod{p^3}$.
\end{rmk}

Let $P_{n,k}$ be the $k$th elementary symmetric polynomial in $n$ variables with values $\{1,2,\dots,n\}$. Since the recurrence relation \cite[equation 14.3]{Gould} for $P_{n,k}$ is \[P_{n,k} = P_{n-1,k}+nP_{n-1,k-1}\] and we have \begin{align} S_{n,n-k} = P_{n,k}/n!, \quad k=0,1,2,\dots,n, \label{form2}
\end{align} we get the recurrence relation \begin{align}S_{n,n-k} = S_{n-1,n-k}+\frac{S_{n-1,n-(k+1)}}{n}. \label{rec} \end{align} So by \eqref{form}, \begin{align} (w_p-1)/p^3 = S_{p,2}/p + pS_{p,4} + p^3S_{p,6} + \cdots + p^{p-6}S_{p,p-3} + p^{p-4}S_{p,p-1}. \label{int} \end{align}
Notice \eqref{fra} and \eqref{rec} imply $S_{p,m} \equiv 0 \pmod{p}$ for $m=2,4,\dots,p-3$. We next seek an explicit formula for the rationals $S_{p,2}, S_{p,4}, \dots, S_{p,p-1}$. 

The \textit{Stirling numbers of the first kind and second kind}, denoted $s(n,k)$ and $S(n,k)$ respectively, are characterized by \[\sum_{k=0}^n s(n,k)x^k = n!\binom{x}{n}, \quad x \in \R,\] \[\sum_{k=0}^n k!\binom{x}{k}S(n,k) = x^n, \quad x \in \R,\] and we have an explicit formula \cite[equation 13.32]{Gould} relating the two sets of numbers, \begin{align} s(n,n-k) = \sum_{j=0}^k(-1)^j\binom{n+j-1}{k+j}\binom{n+k}{k-j}S(j+k,j). \label{form3} \end{align}

\begin{rmk}
Since $S(n,0) = s(n,0) = 0$ for $n \geq 1$, \eqref{form3} and the summations that follow may start at index one. 
\end{rmk}
\noindent
The integers $s(n,k)$ and $P_{n,k}$ are related \cite[equation 14.5]{Gould} by \[P_{n,k} = (-1)^k s(n+1,n+1-k),\] so by \eqref{form2} and \eqref{form3}, the explicit formula for $S_{p,2}, S_{p,4}, \dots, S_{p,p-1}$ is \[S_{p,p-k} = \frac{1}{p!}\sum_{j=0}^k(-1)^{j+1}\binom{p+j}{k+j}\binom{p+1+k}{k-j}S(j+k,j), \quad k=1,3,5,\dots,p-2.\]
Also \begin{align} S_{p,p-k} &= \frac{(p+1+k)!}{p!(p-k)!}\sum_{j=0}^k(-1)^{j+1}\frac{1}{(k+j)!(k-j)!(p+1+j)}S(j+k,j)\nonumber\\
&= \frac{p+1}{(2k)!(p-k)!}\sum_{j=0}^k(-1)^{j+1}\binom{2k}{k+j}\frac{C(p,k)}{p+1+j}S(j+k,j)
\label{form4} \end{align}
where $C(p,k) = (p+2)(p+3)\cdots(p+1+k)$, so \eqref{int} may be written as \begin{align*} \frac{w_p-1}{p^3} &=\frac{p+1}{(2p-4)!(p-1)!}(a_{2p-7}p^{2p-7}+\cdots+a_2p^2+a_1p+a_0)\\
&=\frac{p+1}{(2p-4)!(p-1)!} W(p), \quad a_i \in \Z.
\end{align*}
The properties of some of the coefficients of the \textit{Wolstenholme polynomials} $W(p)$ may easily be determined. In particular $a_{2p-7}$, $a_{2p-8}$, $a_1$, and $a_0$ come from $k=p-2$ in \eqref{form4}, that is from $S_{p,2}$ alone. So in finding the common denominator above, we see $(p-3)!\mid a_0$. We also claim $a_{2p-7} = (2p-5)!!$. By formula \eqref{form3}, \begin{align*}(-1)^ks(n,n-k) &= \sum_{j=0}^k(-1)^{j+k}\binom{n+j-1}{k+j}\binom{n+k}{k-j}S(j+k,j)\\
&= \frac{(n+k)!}{(2k)!(n-1-k)!}\sum_{j=0}^k(-1)^{j+k}\binom{2k}{k+j}\frac{S(j+k,j)}{n+j}\\
&= \frac{1}{(2k)!}\sum_{j=0}^k(-1)^{j+k}\binom{2k}{k+j}\frac{D(n,k)}{n+j}S(j+k,j)
\end{align*}
with $D(n,k) = (n-k)(n-k+1)\cdots (n+k)$. Also from \cite[proposition 1.1]{Gessel}, for each fixed integer $k\geq 0$, $(-1)^ks(n,n-k)$ is a polynomial in $n$ over the rationals with leading coefficient $(2k-1)!!/(2k)!$. Since $D(n,k)$ is a monic polynomial in $n$, \begin{align} (2k-1)!! = \sum_{j=0}^k(-1)^{j+k}\binom{2k}{k+j}S(j+k,j), \label{ident}\end{align} and since $C(p,k)$ is a monic polynomial in $p$ for fixed $k$, \eqref{form4} and \eqref{ident} imply $a_{2p-7} = (2p-5)!!$. Data also suggests $a_{p-4}$ is the largest coefficient and the sign of the coefficients between $a_{2p-7}$ and $a_{p-4}$ alternate.

Let $p$ and $q$ be primes such that $p < q$ and $q\mid(w_p-1)$. Since $q\not\vert ~ w_p$, the proof of \cite[proposition 5, part 4]{Helou} implies $q > 2p$. Hence a prime $q > p$ divides $(w_p-1)/p^3$ if and only if $q\mid W(p)$. This makes the Wolstenholme polynomials useful in situations like conjecture \ref{conj}. For a prime $p$, the Taylor series expansion of $W(p)$ about $n \in \Z$ is \[W(p) = W(n) + W'(n)(p-n) +  \cdots + \frac{W^{(2p-6)}(n)}{(2p-6)!}(p-n)^{2p-6},\] so $(p-n)\mid W(p)$ if and only if $(p-n)\mid W(n)$. Also when $(p-n)\mid W'(n)$, $(p-n)^2\mid W(p)$ if and only if $(p-n)^2\mid W(n)$. If $(p-n)\not\vert \,\, W'(n)$ and $p-n$ is prime, then by Hensel's lemma \cite[theorem 3.4.1]{Gouvea}, $(p-n)^2\mid W(p)$ if and only if $p$ is in the congruence class $s$ mod $(p-n)^2$ where $s \in \Z$ such that $(p-n)^2\mid W(s)$ and $(p-n)\mid(s-n)$. In the case where $p-n$ is a prime larger than $p$ with $(p-n)\mid W(n)$, data collected for primes $p \leq 200$ shows $(p-n)\not\vert \,\, W'(n)$. Therefore showing this data trend continues for all primes would be an important step towards proving conjecture \ref{conj}.

\section*{Acknowledgements}

All data referenced in this paper was obtained using Mathematica 11 and the author thanks the Institute of Mathematics, Academia Sinica for providing this tool.

\bibliographystyle{elsarticle-harv}

\bibliography{biblio.bib}

\end{document}